\documentclass[reqno]{amsart}
\usepackage{hyperref}
\usepackage{amsmath}
\usepackage{amssymb}
\usepackage{amscd}
\usepackage{graphics}
\usepackage{latexsym}
\usepackage{stmaryrd}

\theoremstyle{plain}
\newtheorem{theorem}{Theorem}[section]
\newtheorem{thm}[theorem]{Theorem}

\newtheorem{lem}[theorem]{Lemma}
\newtheorem{prop}[theorem]{Proposition}

\theoremstyle{definition}
\newtheorem{defn}[theorem]{Definition}

\theoremstyle{remark}

\newcommand{\PP}{\mathbb{P}}

\newcommand{\mc}{\mathcal}

\newcommand{\OO}{\mc{O}}

\newcommand{\Kgnb}[1]{\overline{\mathcal{M}}_{#1}}

\newsavebox{\sembox}
\newlength{\semwidth}
\newlength{\boxwidth}

\newsavebox{\semrbox}
\newlength{\semrwidth}
\newlength{\boxrwidth}

\title{Higher Fano manifolds and rational surfaces} 
\author[de Jong]{A. J. de Jong}
\author[Starr]{Jason Starr} 
\begin{document}

%%%%%%%%%%%%%%%%%%%%%%%%%%%%%%%%%%%%%%%%%%%%%%%%%%%%%%%%%%%%%%%%%%%%
%%
%% Abstract
%%
%%%%%%%%%%%%%%%%%%%%%%%%%%%%%%%%%%%%%%%%%%%%%%%%%%%%%%%%%%%%%%%%%%%%

\begin{abstract}
Let $X$ be a Fano manifold of pseudo-index $\geq 3$ such that
$c_1(X)^2-2c_2(X)$ is nef.  Irreducibility of some spaces of rational
curves on $X$ (in fact, a weaker hypothesis) implies a general point
of $X$ is contained in a rational surface.
\end{abstract}

%%%%%%%%%%%%%%%%%%%%%%%%%%%%%%%%%%%%%%%%%%%%%%%%%%%%%%%%%%%%%%%%%%%%%%
%%
%% Body
%%
%%%%%%%%%%%%%%%%%%%%%%%%%%%%%%%%%%%%%%%%%%%%%%%%%%%%%%%%%%%%%%%%%%%%%%

\maketitle

%%%%%%%%%%%%%%%%%%%%%%%%%%%%%%%%%%%%%%%%%%%%%%%%%%%%%%%%%%%%%%%%%%%%%%
%%
%% Section 1.  Introduction
%%
%%%%%%%%%%%%%%%%%%%%%%%%%%%%%%%%%%%%%%%%%%%%%%%%%%%%%%%%%%%%%%%%%%%%%%

\section{Introduction} \label{sec-intro}

One consequence of the bend-and-break lemma is uniruledness of Fano
manifolds, ~\cite{Miyaoka-Mori86}.  In fact, in characteristic $0$,
Fano manifolds are rationally connected, ~\cite{KMM92c}, ~\cite{Ca}.
We prove an analogous theorem with rational curves replaced by
rational surfaces for Fano manifolds satisfying positivity of the
second graded piece of the Chern character.

\begin{defn} \label{defn-1}
A Fano manifold is \emph{$2$-Fano} if $\text{ch}_2(T_X)$ is nef, where
$\text{ch}_2(T_X)$ is the second graded piece of the Chern character,
$\frac{1}{2}(c_1(T_X)^2-2c_2(T_X))$.  In other words,
$\text{deg}(\text{ch}_2(T_X)|_S)$ is nonnegative for every surface $S$
in $X$.
\end{defn}

Let $\mc{M}$ be a positive-dimensional, irreducible component of the
Artin stack $\Kgnb{0,0}(X)$ of genus $0$ stable maps to $X$ whose
general point of $\mc{M}$ parametrizes a stable map with irreducible
domain.  Denote by $M$ the coarse moduli space of $M$.  Denote by
$\Delta$ the locally principal closed substack of $\Kgnb{0,0}(X)$
parametrizing stable maps with reducible domain.  The closed substack
$\mc{M}\cap \Delta$ is a Cartier divisor.  The question we consider is
uniruledness of $M$.

\begin{thm} \label{thm-main}
If $X$ is $2$-Fano, every point of $M$ parametrizing a free curve and
contained in a proper curve in $M - M\cap \Delta$ is contained in a
rational curve in $M$.

If a general point of $M$ parametrizes a birational, free curve and is
contained in a proper curve in $M-M\cap \Delta$, then a general point
of $X$ is contained in a rational surface.
\end{thm}

The proof uses the bend-and-break approach of ~\cite{Miyaoka-Mori86}.
Given a general curve $C$ in $M$, we need to bound the dimension of
$\text{Hom}(C,M)$ from below.  Although $M$ and $\mc{M}$ may be very
singular, the deformation theory of stable maps nonetheless gives a
useful lower bound.  This uses Grothendieck-Riemann-Roch computations
from ~\cite{dJS3}.  Unfortunately, the formula has a negative term
coming from intersection points of $C$ and $\Delta$.  This is the
reason for the hypothesis that $M-M\cap \Delta$ contains a proper
curve.  Luckily, there are nice sufficient conditions for $M-M\cap
\Delta$ to contain many proper curves.

\begin{prop} \label{prop-cont}
If the pseudo-index of $X$ is $\geq 3$ and every irreducible component
of $\mc{M}\cap \Delta$ is an irreducible component of $\Delta$, then
$M - M\cap \Delta$ is a union of proper curves.
\end{prop}

This uses a contraction of the locally principal closed subspace
$\Delta$ in $\Kgnb{0,0}(X)$ discovered in ~\cite{CHS1} and
independently by Adam Parker ~\cite{Parker}.

Section ~\ref{sec-ex} gives some examples of 2-Fano manifolds and
makes some observations about classification.  Section
~\ref{sec-sharp} shows Theorem ~\ref{thm-main} is sharp in 2 ways.
First, there are Fano manifolds that are not $2$-Fano where the
components $\mc{M}$ are not uniruled.  Second, there are $2$-Fano
manifolds where the components $\mc{M}$ are uniruled but not
rationally connected.  Finally Section ~\ref{sec-spec} speculates on
sufficient conditions for the components $\mc{M}$ to be rationally
connected.

\section{Proof of the theorem} \label{sec-pf}

For every point $x$, denote by $\text{Hom}(\PP^1,X,0\mapsto
x)_{\text{nc}}$ the open subscheme of $\text{Hom}(\PP^1,X,0\mapsto x)$
parametrizing nonconstant morphisms.

\begin{lem} \label{lem-1}
The dimension of every irreducible component of
$\text{Hom}(\PP^1,X,0\mapsto x)_{\text{nc}}$ is at least as large as
the pseudo-index of $X$.
\end{lem}

\begin{proof}
This follows from ~\cite[Theorem II.1.2, Corollary II.1.6]{K}.  
\end{proof}

\begin{proof}[Proof of Proposition ~\ref{prop-cont}]
Let $f:X\hookrightarrow \PP^N$ be a plurianticanonical embedding.
Denote by $\Kgnb{0,0}(f):\Kgnb{0,0}(X)\rightarrow \Kgnb{0,0}(\PP^N)$
the associated embedding.  Denote by
$\phi:\Kgnb{0,0}(\PP^N)\rightarrow Y$ the contraction of the boundary
constructed in ~\cite{CHS1}.  Denote by $N$ the image of $\mc{M}$ in
$Y$.

Since the restriction of $\phi$ to $\Kgnb{0,0}(\PP^N) - \Delta$ is an
open immersion, the restriction of $\phi\circ \Kgnb{0,0}(f)$ to
$\mc{M}-\Delta$ is and immersion.  Since $\mc{M}_{\text{free}}$ is
dense in $\mc{M}$, $\mc{M}$ has pure dimension equal to the expected
dimension, and $\mc{M}\cap \Delta$ is a Cartier divisor.  Therefore
$\text{dim}(N)$ equals $\text{dim}(M)$ and $\text{dim}(\mc{M}\cap
\Delta)$ equals $\text{dim}(M)-1$.

If $i\leq j$, the restriction of $\phi$ to the boundary divisor
$\Delta_{i,j}$ factors through the projection
$\pi_j:\Delta_{i,j}\rightarrow \Kgnb{0,1}(\PP^N,j)$.  Denote
$\Delta_{i,j} \cap \Kgnb{0,0}(X)$ by $\Delta_{X,i,j}$.  Denote the
restriction of $\pi_j$ by $\pi_{X,j}:\Delta_{X,i,j} \rightarrow
\Kgnb{0,1}(\PP^N,j)$.  By Lemma ~\ref{lem-1}, every irreducible
component of every fiber of $\pi_{X,j}$ has dimension $\geq 1$, i.e.,
the difference of the pseudo-index and
$\text{dim}(\text{Aut}(\PP^1,0))$.  Therefore, for every irreducible
component $\Delta'$ of $\Delta$, the dimension of
$\phi(\Kgnb{0,0}(f)(\Delta'))$ is strictly less than the dimension of
$\Delta'$.  By hypothesis, every irreducible component $\Delta'$ of
$\mc{M}\cap \Delta$ is an irreducible component of $\Delta$.  Since
$\text{dim}(\Delta')$ equals $\text{dim}(\mc{M})-1$, the image of
$\Delta'$ in $N$ has dimension $\leq \text{dim}(N)-2$.

Since every connected component of $Y$ is projective, also $N$ is
projective.  Because $\text{dim}(\text{Image}(\Delta')) \leq
\text{dim}(N)-2$, a general intersection of $N$ with $\text{dim}(N)-1$
hyperplanes containing a point of $N-\text{Image}(\Delta')$ is a
complete curve that does not intersect $\text{Image}(\Delta')$.
Because there are only finitely many irreducible components of
$\mc{M}\cap \Delta$, a general intersection of $N$ with
$\text{dim}(N)-1$ hyperplanes containing a point of
$N-\text{Image}(\mc{M}\cap \Delta)$ is a complete curve that does not
intersect $\text{Image}(\mc{M}\cap \Delta)$.  The inverse image of
this curve in $\mc{M}-\mc{M}\cap \Delta$ is a complete curve
containing a given point of $\mc{M}-\mc{M}\cap \Delta$.
\end{proof}

Let $C$ be a smooth, proper, connected curve and let
$\zeta:C\rightarrow \Kgnb{0,0}(X)-\Delta$ be a nonconstant 1-morphism
whose general point parametrizes a free curve of $(-K_X)$-degree $e$.
Let $B$ be a finite set of closed points of $C$.  Denote by
$(\pi:\Sigma \rightarrow C,F:\Sigma \rightarrow X)$ the associated
family of stable maps.

\begin{lem} \label{lem-3}
The dimension at $[\zeta]$ of $\text{Hom}(C,\Kgnb{0,0}(X),\zeta|_B)$
is at least,
$$
\text{deg}(\text{ch}_2(T_X)|_{F(\Sigma)}) +
\frac{1}{2e}\text{deg}(c_1(T_X)^2|_{F(\Sigma)}) +
(e+\text{dim}(X)-3)(1-g(C)-\#(B)).
$$
\end{lem}

\begin{proof}
Consider the finite morphism $(\pi,g):\Sigma \rightarrow C\times X$.
Denote by $\mc{N}$ the cokernel of the map,
$$
d(\pi,g):T_\Sigma \rightarrow \pi^* T_C\oplus g^*T_X.
$$
By a natural generalization of ~\cite[Theorem I.2.16]{K}, the
dimension of $\text{Hom}(C,\Kgnb{0,0}(X), \zeta|_B)$ at $\zeta$ is at
least,
$$
h^0(\Sigma,\mc{N}) - h^1(\Sigma,\mc{N}).
$$

By the Leray spectral sequence, $h^2(\Sigma,\mc{N})$ equals
$h^1(C,R^1\pi_* \mc{N})$.  Because a general point of $C$ parametrizes
a free curve, the restriction of $\mc{N}$ to a general fiber of $\pi$
is generated by global sections, thus has no higher cohomology.  Thus
$R^1\pi_*\mc{N}$ is a torsion sheaf so that $h^1(C,R^1\pi_* \mc{N})$
is $0$.  Therefore, the lower bound actually equals
$\chi(\Sigma,\mc{N})$.

Finally, by the Grothendieck-Riemann-Roch computations from
~\cite{dJS3}, $\chi(\Sigma,\mc{N})$ equals,
$$
\text{deg}(\text{ch}_2(T_X)|_{F(\Sigma)}) +
\frac{1}{2e}\text{deg}(c_1(T_X)^2|_{F(\Sigma)}) +
(e+\text{dim}(X)-3)(1-g(C)-\#(B)).
$$
\end{proof}

\begin{proof}[Proof of Theorem ~\ref{thm-main}]
Every proper curve in $M-M\cap \Delta$ is the image of a nonconstant
1-morphism $\zeta:C\rightarrow \mc{M}-\mc{M}\cap \Delta$ from a smooth
curve $C$.  The induced morphism $\text{Hom}(C,\mc{M}-\mc{M}\cap
\Delta) \rightarrow \text{Hom}(C,M)$ is finite.  By Lemma~\ref{lem-3},
$\text{dim}(\text{Hom}(C,M;\zeta|_B))$ behaves as if $M$ is smooth
along the image of $\zeta$ and the anticanonical degree of $\zeta(C)$
equals
$$
\text{deg}(\text{ch}_2(T_X)|_{F(\Sigma)}) +
\frac{1}{2e}\text{deg}(c_1(T_X)^2|_{F(\Sigma)}).
$$ 
Because $X$ is $2$-Fano, this degree is positive.  Therefore the usual
bend-and-break argument applies, cf. ~\cite[Theorem II.5.8]{K}.
\end{proof}

\section{Examples of $2$-Fano manifolds} \label{sec-ex}

All the results of this section, and more, are discussed and proved in
the note ~\cite{dJS6}.  There are a two families of $2$-Fano
manifolds.  The first family comes from complete intersections.  Let
$\PP$ be a weighted projective space of dimension $n$.  Let $X\subset
\PP$ be a smooth complete intersection of type $(d_1,\dots,d_r)$.
Then $X$ is Fano if and only if $d_1 +\dots +d_r \leq n$.  It is
$2$-Fano if and only if $d_1^2 + \dots + d_r^2 \leq n$.

The second family comes from Grassmannians.  Let $\mathbb{G}$ be a
Grassmannian $\text{Grass}(k,n)$ of $k$-dimensional subspaces of a
fixed $n$-dimensional vector space.  Without loss of generality,
assume $n\geq 2k$.  This is Fano.  It is $2$-Fano if and only if
either $k=1$, $n=2k$ or $n=2k+1$.

There are two operations for producing new $2$-Fano manifolds.  First,
if $X$ and $Y$ are each $2$-Fano, then the product $X\times Y$ is
$2$-Fano.  The second operation is more interesting.  Let $X$ be a
smooth Fano manifold and let $L$ be a nef invertible sheaf.  The
$\PP^1$-bundle $\PP(\OO_X\oplus L^\vee)$ is Fano if and only if
$c_1(T_X)-c_1(L)$ is ample.  Assuming it is Fano, it is $2$-Fano if
and only if $\text{ch}_2(T_X) + \frac{1}{2}c_1(L)^2$ is nef.  Notice,
it is not necessary that $\text{ch}_2(T_X)$ is nef, i.e., $X$ need not
be $2$-Fano.

There are other operations on Fano manifolds.  It is reasonable to ask
which of these produce $2$-Fano manifolds.  For instance, a projective
bundle $\PP(E)$ of fiber dimension $\geq 2$ over a Fano manifold is
also Fano if $E$ satisfies a weak version of stability.  However, if
$\PP(E)$ is $2$-Fano then the pullback of $E$ to every curve is a
semistable bundle.  If $X$ is $\PP^n$, for instance, this implies
$\PP(E)$ is simply $\PP^m\times \PP^n$.  This, and other examples,
suggest the following principle: an operation on Fano manifolds
produces a $2$-Fano manifold only if some vector bundle associated to
the operation is semistable.

\section{The theorem is sharp} \label{sec-sharp}

The theorem is sharp in 2 ways.  First, let $X$ be a general cubic
hypersurface in $\PP^5$.  This is Fano, but it is not $2$-Fano.  By
the main theorem of ~\cite{dJS2}, there are infinitely many
non-uniruled irreducible components $\mc{M}$ of $\Kgnb{0,0}(X)$
satisfying the hypotheses of Theorem ~\ref{thm-main}.

Second, let $Y$ be the $\PP^1$-bundle over $X$, $Y = \PP(\OO_X \oplus
\OO_{\PP^5}(-2)|_X)$.  By the construction in the last section, $Y$ is
$2$-Fano.  Associated to the projection $\pi:Y\rightarrow X$, there is
a 1-morphism $\Kgnb{0,0}(\pi):\Kgnb{0,0}(Y)\rightarrow \Kgnb{0,0}(X)$.
For an irreducible component $\mc{N}$ of $\Kgnb{0,0}(Y)$ containing a
free curve, it is easy to prove the boundary $\mc{N}\cap \Delta$ is
contracted.  (However it is not true that every component of
$\mc{N}\cap \Delta$ is a component of $\Delta$.)  Thus Theorem
~\ref{thm-main} implies $N$ is uniruled.  In fact, the restriction of
$\Kgnb{0,0}(\pi)$ to $\mc{N}$ is birational to a projective bundle
over the image component $\mc{M}$ of $\Kgnb{0,0}(X)$.  Choosing
$\mc{N}$ appropriately, $\mc{M}$ is one of the infinitely many
non-uniruled irreducible components of $\Kgnb{0,0}(X)$.  Therefore $N$
is not rationally connected, and the MRC quotient of $N$ is precisely
$M$.

\section{Speculation} \label{sec-spec}

For the counterexample $Y$ in the previous section, $\text{ch}_2(T_Y)$
is nef.  But it is not ``positive''.  It has intersection number $0$
with the surface $\pi^{-1}(B)$ for every curve $B$ in $X$.  If $X$ is
a Fano manifold such that $\text{ch}_2(T_X)$ has positive intersection
number with every surface, is $\mc{M}$ rationally connected?  We know
no counterexample.

\bibliography{my}
\bibliographystyle{alpha}

\end{document}